\documentclass{article}

\usepackage[a4paper]{geometry}
\usepackage[indent]{parskip}
\usepackage{etoolbox}
\usepackage{ifdraft}

\usepackage[dvipsnames]{xcolor}
\definecolor{Basil}{HTML}{32612D}
\definecolor{Berry}{HTML}{7A1712}

\usepackage[author={}]{fixme}
\fxsetup{marginface=\normalfont\scriptsize\sffamily\color{gray}}

\usepackage[mathlines]{lineno}
\ifdraft{\linenumbers}{}

\usepackage[american]{babel}
\usepackage{csquotes}

\input{ushyphex.tex}

\usepackage{lmodern}
\usepackage[
    sf,
    bf,
    nobottomtitles*,
]{titlesec}
\usepackage{abstract}
\titleformat{\section}{\normalfont\sffamily\large\bfseries}{\thesection}{1em}{}
\titleformat{\subsection}{\normalfont\sffamily\normalsize\bfseries}{\thesubsection}{1em}{}

\usepackage[
    useregional,
    showzone=false,
]{datetime2}

\usepackage[nodisplayskipstretch]{setspace}
\usepackage{microtype}
\usepackage[prevent-all]{widows-and-orphans}
\usepackage{titling}
\pretitle{\begin{center}\Large}
\posttitle{\par\end{center}\vskip 1em}
\preauthor{\begin{center}\large \lineskip 0.5em\begin{tabular}[t]{c}}
\postauthor{\end{tabular}\par\end{center}\vskip 1em}
\newcommand{\keywords}[1]{
    \begingroup
    \def\and{ ; }
    \hypersetup{pdfkeywords={#1}}
    \def\and{\ifhmode\unskip\nobreak\fi\ $\cdot$ }
    \paragraph*{Keywords} #1
    \endgroup
}
\newcommand{\msc}[1]{
    \begingroup
    \def\and{ ; }
    \hypersetup{pdfsubject={#1}}
    \def\and{\ifhmode\unskip\nobreak\fi\ $\cdot$ }
    \paragraph*{Mathematics Subject Classification (2020)} #1
    \endgroup
}

\usepackage{fancyhdr}
\fancypagestyle{plain}{
    \fancyhf{}
    
    \fancyfoot[C]{\small\thepage}
}
\fancypagestyle{main}{
    \fancyhf{}
    
    \fancyhead[L]{\small\ifnumodd{\value{page}}{}{\thepage}}
    \fancyhead[R]{\small\ifnumodd{\value{page}}{\thepage}{}}
    \fancyhead[C]{\small\textsc{\ifnumodd{\value{page}}{\headtitle}{\headauthor}}}
}
\usepackage{graphicx}
\usepackage{booktabs}
\usepackage{multirow}
\usepackage{caption}
\usepackage{subcaption}
\graphicspath{{figures/}}
\usepackage{enumitem}
\setlist[itemize]{noitemsep}
\setlist[enumerate]{noitemsep}

\usepackage{siunitx}
\usepackage[
    ruled,
    lined,
    longend,
    linesnumbered,
]{algorithm2e}
\SetKwInput{KwNote}{Note}
\AtBeginEnvironment{algorithm}{\DontPrintSemicolon}
\makeatletter

\newenvironment{algomathdisplay*}{\begin{equation*}}{\@endalgocfline\end{equation*}\ifthenelse{\boolean{algocf@linesnumbered}}{\vspace{-\baselineskip}}{\relax}\;}

\makeatother

\usepackage[final]{listings}
\usepackage{lstautogobble}
\usepackage{tikz}
\usepackage{pgfplots}
\usepackage{pgfplotstable}
\pgfplotsset{compat=1.18}

\usepackage[
    style=ext-numeric-comp,
    articlein=false,
    innamebeforetitle=true,
    sorting=nyt,
    sortcase=false,
    sortcites,
    maxbibnames=99,
    date=year,
    urldate=long,
]{biblatex}
\DeclareNameAlias{author}{family-given}
\renewbibmacro*{edition}{\iffieldundef{edition}{}{\printfield{edition} edition}}

\usepackage{mathtools}
\usepackage{amssymb}
\usepackage{amsthm}
\usepackage{cases}
\usepackage{empheq}
\usepackage[bb=dsfontserif]{mathalpha}
\usepackage{etoolbox}
\usepackage{ifthen}
\usepackage{xargs}
\usepackage{relsize}
\theoremstyle{definition}

\theoremstyle{plain}

\theoremstyle{remark}

\mathcode`@="8000{\catcode`\@=\active\gdef@{\mkern1mu}}

\DeclareMathOperator{\rank}{rank}
\DeclareMathOperator{\sgn}{sgn}

\newcommand*{\R}{\mathbb{R}}

\newcommand*{\T}{\mathsf{T}}
\newcommand*{\abs}[2][]{#1\lvert#2#1\rvert}

\newcommandx*{\norm}[3][1={},2={}]{#1\lVert#3#1\rVert\ifblank{#2}{}{_{\ifthenelse{\equal{#2}{fro}}{\mathsf{F}}{\ifthenelse{\equal{#2}{nuc}}{\ast}{#2}}}}}

\newcommand*{\set}[2][]{#1\{#2#1\}}

\newcommand{\NN}{\mathrm{N}}
\newcommand{\base}{{\text{b}}}
\newcommand{\drop}{{\text{d}}}
\newcommand{\geo}{{\text{g}}}
\newcommand{\st}{\text{s.t.}}
\newcommand{\trust}{{\text{t}}}
\newcommand{\aeq}{A_{\scriptscriptstyle\text{E}}}
\newcommand{\aub}{A_{\scriptscriptstyle\text{I}}}
\newcommand{\beq}{b_{\scriptscriptstyle\text{E}}}
\newcommand{\bub}{b_{\scriptscriptstyle\text{I}}}
\newcommand{\con}[1][i]{c\ifthenelse{\equal{#1}{}}{}{_{#1}}}
\newcommand{\ceq}{\con[\scriptscriptstyle{\text{E}}]}
\newcommandx{\conm}[2][1=k,2=i]{\con[#1, #2]}
\newcommand{\cub}{\con[\scriptscriptstyle\text{I}]}
\newcommand{\eqdef}{\stackrel{\mathsmaller{\mathsf{def}}}{=}}
\newcommand{\frob}{\mathsf{F}}
\newcommand{\fsetm}[1][k]{\Omega_{#1}}
\newcommand{\fset}{\Omega}
\newcommand{\func}{\mathcal{F}}
\newcommand{\iter}[1][k]{x\ifthenelse{\equal{#1}{}}{}{_{#1}}}
\newcommand{\objm}[1][k]{\obj\ifthenelse{\equal{#1}{}}{}{_{#1}}}
\newcommand{\obj}{f}
\newcommand{\qspace}[1][n]{\mathcal{Q}\ifthenelse{\equal{#1}{}}{}{_{#1}}}
\newcommand{\radalt}[1][k]{\widetilde{\Delta}\ifthenelse{\equal{#1}{}}{}{_{#1}}}
\newcommand{\radlb}[1][k]{\rho\ifthenelse{\equal{#1}{}}{}{_{#1}}}
\newcommand{\rad}[1][k]{\Delta\ifthenelse{\equal{#1}{}}{}{_{#1}}}

\newcommand{\xl}{l}
\newcommand{\xpt}[1][k]{\mathcal{Y}\ifthenelse{\equal{#1}{}}{}{_{#1}}}
\newcommand{\xu}{u}

\usepackage[
    pdfusetitle,
    hyperfootnotes=false,
]{hyperref}
\usepackage{url}
\hypersetup{
    final,
    colorlinks,
    linkcolor=Berry,
    citecolor=Basil,
    urlcolor=RedViolet,
}
\newcommand*{\email}[1]{\href{mailto:#1}{\texttt{#1}}}
\newcommand*{\orcid}[1]{\href{https://orcid.org/#1}{\texttt{#1}}}

\usepackage[
    capitalize,
    noabbrev,
]{cleveref}

\crefname{equation}{}{}
\Crefname{equation}{}{}
\crefname{lstlisting}{Listing}{Listings}
\Crefname{lstlisting}{Listing}{Listings}

\usepackage[acronym]{glossaries-extra}
\setabbreviationstyle[acronym]{long-short}
\glsdisablehyper
\newacronym{dfo}{DFO}{derivative-free optimization}
\newacronym{mdo}{MDO}{multidisciplinary design optimization}
\newacronym{psb}{PSB}{Powell's symmetric Broyden}
\newacronym{rbf}{RBF}{radial basis function}
\newacronym{rs}{RS}{Random Search}
\newacronym{tcg}{TCG}{truncated conjugate gradient}
\newacronym{tpe}{TPE}{Tree-Structured Parzen Estimator}
\newglossaryentry{bobyqa}{
    name=BOBYQA,
    description={Bound Optimization BY Quadratic Approximation},
}
\newglossaryentry{bfgs}{
    name=BFGS,
    description={Broyden-Fletcher-Goldfarb-Shanno},
}
\newglossaryentry{cg}{
    name=CG,
    description={Conjugate Gradient},
}
\newglossaryentry{cobyla}{
    name=COBYLA,
    description={Constrained Optimization BY Linear Approximations},
}
\newglossaryentry{cobyqa}{
    name=COBYQA,
    description={Constrained Optimization BY Quadratic Approximations},
}
\newglossaryentry{lincoa}{
    name=LINCOA,
    description={LINearly Constrained Optimization Algorithm},
}
\newglossaryentry{newuoa}{
    name=NEWUOA,
    description={NEW Unconstrained Optimization Algorithm},
}
\newglossaryentry{pdfo}{
    name=PDFO,
    description={Powell's Derivative-Free Optimization solvers},
}
\newglossaryentry{prima}{
    name=PRIMA,
    description={Reference Implementation for Powell's methods with Modernization and Amelioration},
}
\newglossaryentry{uobyqa}{
    name=UOBYQA,
    description={Unconstrained Optimization BY Quadratic Approximation},
}

\NewCommandCopy{\oldtitle}{\title}
\newcommand{\headtitle}{}
\renewcommand{\title}[2][]{
    \oldtitle{#2}
    \renewcommand{\headtitle}{\ifblank{#1}{#2}{#1}}
}

\NewCommandCopy{\oldauthor}{\author}
\newcommand{\headauthor}{}
\renewcommand{\author}[2][]{
    \oldauthor{#2}
    \renewcommand{\headauthor}{%
        \renewcommand{\thanks}[1]{}%
        \renewcommand{\thanksgap}[1]{}%
        \renewcommand{\thanksmark}[1]{}%
        \renewcommand{\and}{\unskip, }%
        \ifblank{#1}{#2}{#1}%
    }
}

\title[PDFO: Powell's Derivative-Free Optimization solvers]{PDFO: A Cross-Platform Package for Powell's Derivative-Free Optimization Solvers}
\author[T. M. Ragonneau \and Z. Zhang]{%
    Tom M. Ragonneau%
    \thanks{Department of Applied Mathematics, The Hong Kong Polytechnic University, Hong Kong, China.}%
    \thanksgap{1ex}%
    \thanks{Email: \email{tom.ragonneau@gmail.com}; ORCID: \orcid{0000-0003-2717-2876}.}
    \and
    Zaikun Zhang%
    \texorpdfstring{\thanksmark{1}}{}%
    \thanksgap{1ex}%
    \thanks{Email: \email{zaikunzhang@gmail.com}; ORCID: \orcid{0000-0001-8934-8190}.}
}
\date{\ifdraft{\DTMnow}{\today}}

\begin{document}

\maketitle

\begin{abstract}
    The late Professor M. J. D. Powell devised five trust-region methods for derivative-free optimization, namely \gls{cobyla}, \gls{uobyqa}, \gls{newuoa}, \gls{bobyqa}, and \gls{lincoa}.
    He carefully implemented them into publicly available solvers, renowned for their robustness and efficiency.
    However, the solvers were implemented in Fortran 77 and hence may not be easily accessible to some users.
    We introduce the \gls{pdfo} package, which provides user-friendly Python and MATLAB interfaces to Powell's code.
    With \gls{pdfo}, users of such languages can call Powell's Fortran solvers easily without dealing with the Fortran code.
    Moreover, \gls{pdfo} includes bug fixes and improvements, which are particularly important for handling problems that suffer from ill-conditioning or failures of function evaluations.
    In addition to the \gls{pdfo} package, we provide an overview of Powell's methods, sketching them from a uniform perspective, summarizing their main features, and highlighting the similarities and interconnections among them.
    We also present experiments on \gls{pdfo} to demonstrate its stability under noise, tolerance of failures in function evaluations, and potential to solve certain hyperparameter optimization problems.
\end{abstract}

\msc{65K05 \and 90C30 \and 90C56 \and 90-04}

\section{Introduction}

Most optimization algorithms rely on classical or generalized derivative information of the objective and constraint functions.
However, in many applications, such information is not available.
This is the case, for example, if the objective function does not have an explicit formulation but can only be evaluated through complex simulations or experiments.
Such problems motivate the development of optimization algorithms that use only function values but not derivatives, also known as \gls{dfo} algorithms.

Powell devised five algorithms to tackle unconstrained and constrained problems without using derivatives, namely~\gls{cobyla}~\cite{Powell_1994}, \gls{uobyqa}~\cite{Powell_2002}, \gls{newuoa}~\cite{Powell_2006}, \gls{bobyqa}~\cite{Powell_2009}, and \gls{lincoa}.
Not only did he propose these algorithms but he also implemented them into publicly available solvers, paying great attention to the stability and complexity of their numerical linear algebra computations.
Renowned for their robustness and efficiency, these solvers are used in a wide spectrum of applications, for instance, aeronautical engineering~\cite{Gallard_Etal_2018}, astronomy~\cite{Mamon_Biviano_Boue_2013}, computer vision~\cite{Izadinia_Shan_Seitz_2017}, robotics~\cite{Mombaur_Truong_Laumond_2010}, and statistics~\cite{Bates_Etal_2015}.

However, Powell implemented the solvers in Fortran 77, an old-fashioned language that damps the enthusiasm of many users to exploit these solvers in their projects.
There has been a continued demand from both researchers and practitioners for the availability of Powell's solvers in more user-friendly languages such as Python and MATLAB.

Responding to such a demand, this paper presents a package named \gls{pdfo}, an acronym for ``\glsdesc{pdfo}.''
\Gls{pdfo} interfaces Powell's Fortran solvers with other languages, enabling users of such languages to call Powell's solvers without dealing with the Fortran code.
For each supported language, \gls{pdfo} provides a simple function that can invoke one of Powell's solvers according to the user's request (if any) or according to the type of the problem to solve.
The current release (Version~2.2.0) of \gls{pdfo} supports Python and MATLAB, with more languages to be covered in the future.
The signature of the Python function is consistent with the \texttt{minimize} function of the SciPy optimization library, and that of the MATLAB function is consistent with the \texttt{fmincon} function of the MATLAB Optimization Toolbox.
\Gls{pdfo} is cross-platform, available on Linux, macOS, and Windows at
\begin{center}
    \url{https://www.pdfo.net} \quad and \quad \url{https://github.com/pdfo/pdfo}\,,
\end{center}
with the DOI \href{https://doi.org/10.5281/zenodo.3887568}{\texttt{10.5281/zenodo.3887568}}.
It has been downloaded more than \num{120000} times as of June 2024, mirror downloads excluded.
Moreover, it is one of the optimization engines in GEMSEO~\cite{Gallard_Etal_2018},\footnote{\url{https://gemseo.readthedocs.io}\,.} an industrial software package for \gls{mdo}.

\Gls{pdfo} is not the first attempt to facilitate the usage of Powell's solvers in languages other than Fortran.
Various efforts have been made in this direction.
Py-BOBYQA~\cite{Cartis_Etal_2019,Cartis_Roberts_Sheridan-Methven_2022} provides a Python implementation of \gls{bobyqa} (although it is not meant to b
e a faithful re-implementation of Powell's version); NLopt includes multi-language interfaces for \gls{cobyla}, \gls{newuoa}, and \gls{bobyqa};\footnote{\url{https://github.com/stevengj/nlopt}\,.} minqa wraps \gls{uobyqa}, \gls{newuoa}, and \gls{bobyqa} in R;\footnote{\url{https://cran.r-project.org/package=minqa}\,.} SciPy makes \gls{cobyla} available in Python under its optimization library.\footnote{\url{https://docs.scipy.org/doc/scipy/reference/optimize.minimize-cobyla.html}\,.}
However, \gls{pdfo} has several features that distinguish it from others.
\begin{enumerate}
    \item \emph{Comprehensiveness.}
        To the best of our knowledge, \gls{pdfo} is the only package that provides all of \gls{cobyla}, \gls{uobyqa}, \gls{newuoa}, \gls{bobyqa}, and \gls{lincoa} with a \emph{uniform interface}.
    \item \emph{Solver selection.}
        \Gls{pdfo} can automatically select a solver for a given problem.
        The selection takes into account the performance of the solvers on the CUTEst~\cite{Gould_Orban_Toint_2015} problem set.
    \item \emph{Problem preprocessing.}
        \Gls{pdfo} preprocesses the inputs to simplify the problem and reformulate it to meet the requirements of Powell's solvers.
    \item \emph{Code patching.}
        \Gls{pdfo} patches several bugs in the Fortran code.
        Such bugs can lead to serious problems such as infinite cycling or memory errors.
    \item \emph{Fault tolerance.}
        \Gls{pdfo} tolerates failures of function evaluations.
        In case of such failures, \gls{pdfo} will not exit but try to progress.
    \item \emph{Additional options.}
        \gls{pdfo} includes options for the user to control the solvers in some manners that are useful in practice.
        For example, the user can request \gls{pdfo} to scale the problem according to bound constraints on the variables before solving.
\end{enumerate}

In addition to the \gls{pdfo} package, this paper also provides an overview of Powell's \gls{dfo} methods.
We will not repeat Powell's description of these methods but summarize them from a uniform viewpoint, aiming at easing the understanding of Powell's methods and paving the way for further development based on them.

The analysis of Powell's \gls{dfo} methods is not within the scope of this paper.
Under some assumptions, adapted versions of Powell's algorithms may be covered by existing theory for trust-region \gls{dfo} methods~\cite[Chapter~10]{Conn_Scheinberg_Vicente_2009b}.
However, it will still be interesting to pursue a tailored theory for Powell's algorithms.

The remaining part of this paper is organized as follows.
Section~\ref{sec:dfo} briefly reviews \gls{dfo} methods in order to provide the context of \gls{pdfo}.
We then present an overview of Powell's \gls{dfo} methods in Section~\ref{sec:powell}, including a sketch of the algorithms and a summary of their main features.
A detailed exposition of~\gls{pdfo} is given in Section~\ref{sec:pdfo}, highlighting its solver selection, problem preprocessing, bug fixes, and handling of function evaluation failures.
Section~\ref{sec:numerical} presents some experiments on~\gls{pdfo}, demonstrating its stability under noise, tolerance of function evaluation failures, and potential in hyperparameter optimization.
We conclude the paper with some remarks in Section~\ref{sec:conclude}.

\section{A brief review of \texorpdfstring{\gls{dfo}}{DFO} methods}
\label{sec:dfo}

Consider a nonlinear optimization problem
\begin{equation}
    \label{eq:nlc}
    \min_{x \in \fset} ~ \obj(x),
\end{equation}
where~$\obj \colon \R^n \to \R$ is the objective function and~$\fset \subseteq \R^n$ represents the feasible region.
As summarized in~\cite{Conn_Scheinberg_Vicente_2009b}, two strategies have been developed to tackle problem~\eqref{eq:nlc} without using derivatives, which we will introduce in the following.

The first strategy, known as direct search,\footnote{In some early papers (e.g.,~\cite{Powell_1994,Powell_1998}) Powell and many other authors used ``direct search'' to mean what is known as ``\glsxtrlong{dfo}'' today.
Powell rarely used the word ``derivative-free optimization.''
The only exceptions known to us are his last paper~\cite{Powell_2015} and his distinguished lecture titled ``A parsimonious way of constructing quadratic models from values of the objective function in derivative-free optimization'' at the National Center for Mathematics and Interdisciplinary Sciences, Beijing on November 4, 2011~\cite{Buhmann_Fletcher_Iserles_Toint_2018}.} samples the objective function~$\obj$ and chooses iterates by simple comparisons of function values, examples including the Nelder-Mead algorithm~\cite{Nelder_Mead_1965}, the MADS methods~\cite{Audet_Dennis_2006,LeDigabel_2011}, and BFO~\cite{Porcelli_Toint_2017,Porcelli_Toint_2022}.
See~\cite{Kolda_Lewis_Torczon_2003},~\cite[Chapters~7 and~8]{Conn_Scheinberg_Vicente_2009b},~\cite[Part~3]{Audet_Hare_2017}, and~\cite[\S~2.1]{Larson_Menickelly_Wild_2019} for more discussions on this paradigm, and we refer to~\cite{Gratton_Etal_2015,Gratton_Etal_2019} for recent developments of randomized methods in this category.

The second strategy approximates the original problem~\eqref{eq:nlc} by relatively simple models and locates the iterates according to these models.
Algorithms applying this strategy are referred to as model-based methods.
They often make use of the models within a trust-region framework~\cite{Conn_Scheinberg_Vicente_2009a} or a line-search framework~\cite{Berahas_Byrd_Nocedal_2019}.
Interpolation and regression are two common ways of establishing the models~\cite{Powell_2001,Conn_Scheinberg_Vicente_2008a,Conn_Scheinberg_Vicente_2008b,Wild_Regis_Shoemaker_2008,Bandeira_Scheinberg_Vicente_2012,Billups_Larson_Graf_2013,Regis_Wild_2017}.
Algorithms using finite-difference gradients can also be regarded as model-based methods because such gradients essentially come from linear (for the forward and backward differences) or quadratic (for the central difference) interpolation of the function under consideration over rather special interpolation sets~\cite[\S~1.4.3]{Ragonneau_2022}.
Most model-based \gls{dfo} methods employ linear or quadratic models, examples including Powell's algorithms~\cite{Powell_1994,Powell_2002,Powell_2006,Powell_2009} in \gls{pdfo}, MNH~\cite{Wild_2008}, DFLS~\cite{Zhang_Conn_Scheinberg_2010}, DFO-TR~\cite{Bandeira_Scheinberg_Vicente_2012}, and DFO-LS~\cite{Cartis_Etal_2019,Hough_Roberts_2022}, but there are also methods exploiting \glspl{rbf}, such as ORBIT~\cite{Wild_Regis_Shoemaker_2008}, CONORBIT~\cite{Regis_Wild_2017}, and BOOSTERS~\cite{Oeuvray_Bierlaire_2009}.

Hybrids between direct search and model-based approaches exist, for instance, Implicit Filtering~\cite[Algorithm~4.7]{Kelley_2011} and MADS with quadratic models~\cite{Conn_LeDigabel_2013}.
Theories of global convergence and convergence rate have been established for both direct search~\cite{Torczon_1997,Kolda_Lewis_Torczon_2003,Vicente_2013,Gratton_Etal_2015,Dodangeh_Vicente_2016} and model-based methods~\cite{Conn_Scheinberg_Toint_1997a,Conn_Scheinberg_Vicente_2009a,Powell_2012,Garmanjani_Judice_Vicente_2016}.
Since the objective and constraint functions in \gls{dfo} problems are commonly expensive to evaluate, the worst-case complexity in terms of function evaluations is a major theoretical aspect of \gls{dfo} algorithms.
Examples of such complexity analysis can be found in~\cite{Vicente_2013,Gratton_Etal_2015,Dodangeh_Vicente_2016,Garmanjani_Judice_Vicente_2016}.
For more extensive discussions on \gls{dfo} methods and theory, see the monographs~\cite{Conn_Scheinberg_Vicente_2009b,Audet_Hare_2017}, the survey papers~\cite{Rios_Sahinidis_2013,Custodio_Scheinberg_Vicente_2017,Larson_Menickelly_Wild_2019}, the recent thesis~\cite{Ragonneau_2022}, and the references therein.

\section{Powell's derivative-free algorithms}
\label{sec:powell}

Powell published in 1964 his first \gls{dfo} algorithm based on conjugate directions~\cite{Powell_1964}.\footnote{According to Google Scholar, this is Powell's second published paper and also the second most cited work.
The earliest and most cited one is his paper on the DFP method~\cite{Fletcher_Powell_1963} co-authored with Fletcher and published in 1963.
DFP is not a DFO algorithm but the first quasi-Newton method.
The least-change property~\cite{Dennis_Schnabel_1979} of quasi-Newton methods is a major motivation for Powell to investigate the least Frobenius norm updating~\cite{Powell_2004b} of quadratic models in DFO, which is the backbone of \gls{newuoa}, \gls{bobyqa}, and \gls{lincoa}.}
His code for this algorithm is contained in the HSL Mathematical Software Library as subroutine \texttt{VA24}.\footnote{\url{https://www.hsl.rl.ac.uk}\,.}
It is not included in \gls{pdfo} because the code is not in the public domain, although open-source implementations are available (see~\cite[Footnote~4]{Conn_Scheinberg_Toint_1997b}).

From the 1990s to the final years of his career, Powell developed five model-based \gls{dfo} algorithms to solve~\eqref{eq:nlc}, namely \gls{cobyla}~\cite{Powell_1994} (for nonlinearly constrained problems), \gls{uobyqa}~\cite{Powell_2002} (for unconstrained problems), \gls{newuoa}~\cite{Powell_2006} (for unconstrained problems), \gls{bobyqa}~\cite{Powell_2009} (for bound-constrained problems), and \gls{lincoa} (for linearly constrained problems).
Moreover, Powell implemented these algorithms into Fortran solvers and made the code publicly available.
They are the cornerstones of \gls{pdfo}.
This section provides an overview of these five algorithms, starting with a sketch in Section~\ref{ssec:sketch} and then presenting more details afterward.

\subsection{A sketch of the algorithms}
\label{ssec:sketch}

Powell's model-based \gls{dfo} algorithms are trust-region methods.
At iteration~$k$, the algorithms construct a linear (for \gls{cobyla}) or quadratic (for the other methods) model~$\objm$ for the objective function~$f$ according to the interpolation condition
\begin{equation}
    \label{eq:itpls}
    \objm(y) = \obj(y), \quad y \in \xpt,
\end{equation}
where~$\xpt \subseteq \R^n$ is a finite interpolation set updated along the iterations.
\Gls{cobyla} models the constraints by linear interpolants on~$\xpt$ as well.
Instead of repeating Powell's description of these algorithms, we outline them in the sequel, emphasizing the trust-region subproblem, the interpolation problem, and the management of the interpolation set.

\subsubsection{The trust-region subproblem}

In all five algorithms, iteration~$k$ places the trust-region center~$\iter$ at the ``best'' point where the objective function and constraints have been evaluated so far.
Such a point is selected according to the objective function or a merit function that takes the constraints into account.
After choosing the trust-region center~$\iter$, with the trust-region model~$\objm$ constructed according to~\eqref{eq:itpls}, a trial point~$\iter^{\trust}$ is then obtained by solving approximately the trust-region subproblem
\begin{equation}
    \label{eq:trsp}
        \min_{x \in \fsetm} ~ \objm(x) \quad \st \quad \norm{x - \iter} \le \rad,
\end{equation}
where~$\rad$ is the trust-region radius, and~$\norm{\cdot}$ is the~$\ell_2$-norm in~$\R^n$.
In this subproblem, the set~$\fsetm \subseteq \R^n$ is a local approximation of the feasible region~$\fset$.
\Gls{cobyla} defines~$\fsetm$ by linear interpolants of the constraint functions over the set~$\xpt$, whereas the other four algorithms take~$\fsetm = \fset$.

\subsubsection{The interpolation problem}
\label{ssec:iptprob}

\paragraph{Fully determined interpolation.}

The interpolation condition~\eqref{eq:itpls} is essentially a linear system.
Given a base point~$y^{\base}\in \R^n$, which may depend on~$k$, a linear model~$\objm$ takes the form of~$\objm(x) = \obj(y^{\base}) + (x - y^{\base})^{\T} \nabla \objm(y^{\base})$, and hence~\eqref{eq:itpls} is equivalent to
\begin{equation}
    \label{eq:litpls}
    \objm(y^{\base}) + (y -y^{\base})^{\T} \nabla \objm(y^{\base}) = \obj(y), ~ y \in \xpt,
\end{equation}
which is a linear system with respect to~$\objm(y^\base) \in \R$ and~$\nabla \objm(y^\base) \in \R^n$, the degrees of freedom being~$n + 1$.
\Gls{cobyla} builds linear models by the system~\eqref{eq:litpls}, with~$\xpt$ being an interpolation set of~$n+1$ points updated along the iterations.
Similarly, if~$\objm$ is a quadratic model, then~\eqref{eq:itpls} is equivalent to
\begin{equation}
    \label{eq:qitpls}
    \objm(y^{\base}) + (y -y^{\base})^{\T} \nabla \objm(y^{\base}) + \frac{1}{2}(y-y^{\base})^{\T} \nabla^2 \objm(y^{\base}) (y-y^{\base}) = \obj(y), ~ y \in \xpt,
\end{equation}
a linear system with unknowns~$\objm(y^\base) \! \in \! \R$, $\nabla \objm(y^\base) \! \in \! \R^n$, and~$\nabla^2 \objm(y^{\base}) \! \in \! \R^{n \times n}$, the degrees of freedom being~$(n + 1)(n + 2) / 2$ due to the symmetry of~$\nabla^2 \objm(y^\base)$.
\Gls{uobyqa} constructs quadratic models by the system~\eqref{eq:qitpls}.
To decide a quadratic model~$\objm$ completely by this system alone, \gls{uobyqa} requires that~$\xpt$ contains~$(n+1)(n+2)/2$ points, and~$f$ should have been evaluated at all these points before the system can be formed.
Even though most of these points will be reused at the subsequent iterations so that the number of function evaluations needed per iteration is tiny (see Section~\ref{ssec:iptset}), we must perform~$(n + 1)(n + 2) / 2$ function evaluations during the very first iteration.
This is impracticable unless~$n$ is small, which motivates the underdetermined quadratic interpolation.

\paragraph{Underdetermined quadratic interpolation.}

In this case, models are established according to the interpolation condition~\eqref{eq:itpls} with~$\abs{\xpt}$ being less than or equal to~$(n + 1)(n + 2) / 2$, the remaining degrees of freedom being taken up by minimizing a certain functional~$\func_k$ to promote the regularity of the quadratic model.
More specifically, this means building~$\objm$ by solving
\begin{equation}
    \label{eq:undqitp}
        \min_{Q \in \qspace} ~ \func_k(Q) \quad \st \quad Q(y) = \obj(y), ~ y \in \xpt,
\end{equation}
where~$\qspace$ is the space of polynomials on~$\R^n$ of degree at most~\num{2}.
\Gls{newuoa}, \gls{bobyqa}, and \gls{lincoa} construct quadratic models in this way, with
\begin{equation}
    \label{eq:leastchange}
    \func_k(Q) = \norm[\big]{\nabla^2 Q - \nabla^2 \objm[k - 1]}_{\frob}^2,
\end{equation}
which is inspired by the least-change property of quasi-Newton updates~\cite{Dennis_Schnabel_1979}, although other functionals are possible (see, e.g.,~\cite{Conn_Toint_1996,Bandeira_Scheinberg_Vicente_2012,Powell_2013,Zhang_2014,Xie_Yuan_2023}).
The first model~$\objm[1]$ is obtained by setting~$\objm[0] = 0$.
Powell~\cite{Powell_2013} referred to his approach as the \emph{symmetric Broyden update} of quadratic models (see also~\cite[\S~3.6]{Zhang_2012} and~\cite[\S~2.4.2]{Ragonneau_2022}).
It can be regarded as a derivative-free version of \gls{psb} quasi-Newton update~\cite{Powell_1970b}, which minimizes the functional~$\func_k$ among all quadratic polynomials that fulfill~$Q(\iter) = \obj(\iter)$, $\nabla Q(\iter) = \nabla \obj(\iter)$, and~$\nabla Q(\iter[k - 1]) = \nabla \obj(\iter[k - 1])$ (see~\cite[Theorem~4.2]{Dennis_Schnabel_1979}), with~$\iter$ and~$\iter[k - 1]$ being the current and the previous iterates, respectively.
The interpolation problem~\eqref{eq:undqitp}--\eqref{eq:leastchange} is a convex quadratic programming problem with respect to the coefficients of the quadratic model.

\paragraph{Solving the interpolation problem.}

Powell's algorithms do not solve the interpolation problems~\eqref{eq:litpls}, \eqref{eq:qitpls}, and~\eqref{eq:undqitp}--\eqref{eq:leastchange} from scratch.
\Gls{cobyla} maintains the inverse of the coefficient matrix for~\eqref{eq:litpls} and updates it along the iterations.
Since each iteration of \gls{cobyla} alters the interpolation set~$\xpt$ by only one point (see Subsection~\ref{ssec:iptset}), the coefficient matrix is modified by a rank-$1$ update, and hence its inverse can be updated according to the Sherman-Morrison-Woodbury formula~\cite{Hager_1989}.
\Gls{uobyqa} does the same for~\eqref{eq:qitpls}, except that~\cite[\S~4]{Powell_2002} describes the update in terms of the Lagrange functions of the interpolation problem~\eqref{eq:qitpls}, the coefficients of a Lagrange function corresponding precisely to a column of the inverse matrix.
For the underdetermined quadratic interpolation~\eqref{eq:undqitp}--\eqref{eq:leastchange}, \gls{newuoa}, \gls{bobyqa}, and \gls{lincoa} maintain and update the inverse of the coefficient matrix for the KKT system of~\eqref{eq:undqitp}--\eqref{eq:leastchange}.
The update is also done by the Sherman-Morrison-Woodbury formula as detailed in~\cite[\S~2]{Powell_2004c}.
In this case, each iteration modifies the coefficient matrix and its inverse by rank-$2$ updates.
In addition, the columns of this inverse matrix readily provide the coefficients of Lagrange
functions that make the interpolation problem~\eqref{eq:undqitp}--\eqref{eq:leastchange} easy to
solve (see~\cite[\S~3]{Powell_2004b}).

\paragraph{The base point.}

The choice of the base point~$y^{\base}$ is also worth mentioning.
\Gls{cobyla} sets~$y^{\base}$ to the center~$x_k$ of the current trust region.
In contrast, the other four algorithms initiate~$y^{\base}$ to the starting point provided by the user and keep it unchanged except for occasionally updating~$y^{\base}$ to~$\iter$, without which the distance~$\norm{y^{\base}-\iter}$ may become unfavorably large for the numerical solution of the interpolation problem.
See~\cite[\S~5]{Powell_2004b} and~\cite[\S~7]{Powell_2006} for more elaboration.

\subsubsection{The interpolation set}
\label{ssec:iptset}

The strategy to update~$\xpt$ is crucial.
It should reuse points from previous iterations, at which the objective and constraint functions have been evaluated.
Meanwhile, it needs to maintain the geometry of the interpolation set so that it is well poised, or equivalently, the interpolation problem is well conditioned~\cite{Conn_Scheinberg_Vicente_2009b}.

At a normal iteration, Powell's methods compute a point~$\iter^{\trust} \in \R^n$ by solving the trust-region subproblem~\eqref{eq:trsp}, and update the interpolation set as
\begin{equation}
    \label{eq:xpt-update-tr}
    \xpt[k + 1] = \big(\xpt \cup \set{\iter^{\trust}}\big) \setminus \set{y_k^{\drop}},
\end{equation}
where~$y_k^{\drop} \in \xpt$ is selected after obtaining~$\iter^{\trust}$, aiming to maintain the well-poisedness of~$\xpt[k + 1]$.
As mentioned, Powell's methods update the inverse of the coefficient matrix for either the interpolation system or the corresponding KKT system by the Sherman-Morrison-Woodbury formula.
To keep the interpolation problem well-conditioned, $y_k^{\drop}$ is chosen
to enlarge the magnitude of the denominator in this formula, which is also the ratio between the determinants of the old and new coefficient matrices.\footnote{
    Suppose that~$W$ is a square matrix and consider~$\widetilde{W} = W + UV^\T$, where $U$ and~$V$ are two matrices of the same size and $UV^\T$ has the same size as $W$.
    Then $\det(\widetilde{W}) = \det(W)\det(I+V^\T W^{-1}U)$, and the Sherman-Morrison-Woodbury formula is $\widetilde{W}^{-1} = W^{-1} -W^{-1}U(I+V^\T W^{-1}U)^{-1} V^\T W^{-1}$, assuming that both~$W$ and~$I+V^\T W^{-1}U$ are nonsingular.
    The number~$\det(I+V^\T W^{-1}U)$ is the only denominator involved in the numerical computation of the formula.
}
In the fully determined interpolation, this denominator is~$\ell_k^{@\drop}(\iter^{\trust})$, where~$\ell_k^{@\drop}$ is the Lagrange function associated with~$\xpt$ corresponding to~$y_k^\drop$ (see equations~(10)--(13) and~\S~2 of~\cite{Powell_2001}).
In the underdetermined case, the denominator is lower bounded by~$[\ell_k^{@\drop}(\iter^{\trust})]^2$ (see equation~(2.12), Lemma~1, and~\S~2 of~\cite{Powell_2004c}, where the denominator is denoted by $\sigma$, and~$\ell_k^{@\drop}(\iter^{\trust})$ by~$\tau$).
However, Powell's methods do not choose the point~$y_k^\drop$ merely according to this denominator, but also take into account its distance to the trust-region center, giving a higher priority to farther points, as we can see in~\cite[Equation~(56)]{Powell_2002} and~\cite[Equations~(7.4)--(7.5)]{Powell_2006}, for example.

An alternative update of the interpolation set takes place when the methods detect that~$\objm$ does not represent~$\obj$ well enough, attempting to improve the geometry of the interpolation set.
In this case, the methods first select a point~$y_k^{\drop} \in \xpt$ to drop from~$\xpt$, and then set
\begin{equation}
    \label{eq:xpt-update-geo}
    \xpt[k + 1] = \big(\xpt \setminus \set{y_k^{\drop}}\big) \cup \set{\iter^{@\geo}},
\end{equation}
where~$\iter^{@\geo} \in \R^n$ is chosen to improve the well-poisedness of~$\xpt[k + 1]$.
In \gls{cobyla}, the choice of~$y_k^{\drop}$ and~$\iter^{@\geo}$ is guided by the fact that the interpolation set forms a simplex in~$\R^n$, trying to keep~$\xpt[k+1]$ away from falling into an~$(n-1)$-dimensional subspace, as is detailed in~\cite[Equations~(15)--(17)]{Powell_1994}.
The other four methods select~$y_k^{\drop}$ from~$\xpt$ by maximizing its distance to the current trust-region center~$\iter$, and then obtain~$\iter^{@\geo}$ by solving
\begin{equation}
    \label{eq:biglag}
        \max_{x \in \fset} ~ \abs{\ell_k^{@\drop}(x)} \quad \st \quad \norm{x - \iter} \le \radalt
\end{equation}
for some~$\radalt \in (0, \rad]$.
The motivation for this problem is again to enlarge the magnitude of the aforementioned denominator in the Sherman-Morrison-Woodbury updating formula: for~\gls{uobyqa}, the denominator is~$\ell_k^{@\drop}(x)$, while for~\gls{newuoa}, \gls{bobyqa}, and \gls{lincoa}, the denominator is lower bounded by~$[\ell_k^{@\drop}(x)]^2$.
In addition, \gls{newuoa} maximizes this denominator directly if~\eqref{eq:biglag} fails to make its magnitude large enough, which rarely happens~\cite[\S~6]{Powell_2006}.

Given the two possible updates~\eqref{eq:xpt-update-tr} and~\eqref{eq:xpt-update-geo} of the interpolation set, it is clear that the number of interpolation points remains constant.
As mentioned earlier, this number is~$n+1$ in \gls{cobyla} and~$(n+1)(n+2)/2$ in \gls{uobyqa}.
\Gls{newuoa}, \gls{bobyqa}, and \gls{lincoa} set it to an integer in~$[n+2, (n+1)(n+2)/2]$, with the default value being~$2n+1$, which is proved optimal in terms of the well-poisedness of the initial interpolation set chosen by Powell for \gls{newuoa}~\cite{Ragonneau_Zhang_2023a}.

\subsection{\glsfmttext{cobyla}}
\label{ssec:cobyla}

Published in 1994, \gls{cobyla} was the first model-based \gls{dfo} solver by Powell.
The solver is named after ``\glsdesc{cobyla}.''
It aims to solve problem~\eqref{eq:nlc} with the feasible region
\begin{equation*}
    \fset \eqdef \set{x \in \R^n : \con(x) \ge 0, ~ i = 1, \dots, m},
\end{equation*}
where~$\con \colon \R^n \to \R$ denotes the~$i$th constraint function for each $i \in \set{1, \dots, m}$.
The same as the objective function, all constraints are assumed to be accessible only through function values.

As mentioned before, iteration~$k$ of~\gls{cobyla} models the objective and the constraint functions with linear interpolants on the interpolation set~$\xpt$ of~$n + 1$ points.
Once the linear models~$\conm$ of~$\con$ are built for~$i \in \set{1, \dots, m}$, the trust-region subproblem~\eqref{eq:trsp} is formed with
\begin{equation}
    \label{eq:cobylarg}
    \fsetm \eqdef \set{x \in \R^n : \conm(x) \ge 0, ~ i=1, \dots, m}.
\end{equation}
This subproblem may not be feasible, as the trust region and the region~\eqref{eq:cobylarg} may not intersect.
\Gls{cobyla} handles the trust-region subproblem in two stages.
In the first stage, it solves
\begin{equation*}
    \min_{x \in \R^n} \max_{1 \le i \le m} ~ [\conm(x)]_{-} \quad \st \quad \norm{x - \iter} \le \rad,
\end{equation*}
where~$[t]_{-} = \max \set{0, -t}$ for any~$t\in \R$.
In doing so, the method attempts to reduce the~$\ell_{\infty}$-violation of the linearized constraints within the trust region.
If the first stage finds a point in the interior of the trust region, then the second stage uses the resultant freedom in $x$ to minimize the linearized objective function~$\objm$ within the trust region subject to no increase in any greatest violation of the linearized constraints.

\Gls{cobyla} assesses the quality of points and updates the trust-region radius according to an~$\ell_\infty$-merit function and a reduction ratio based on it (see~\cite[Equations~(5),~(9), and~(10)]{Powell_1994}).
It never increases the trust-region radius and reduces it if the geometry of~$\xpt$ is acceptable but the trust-region trial point~$\iter^\trust$ is too close to~$\iter$ or does not render a big enough reduction ratio~\cite[Equation~(11)]{Powell_1994}.

\subsection{\glsfmttext{uobyqa}}
\label{ssec:uobyqa}

In 2002, Powell published \gls{uobyqa}~\cite{Powell_2002}, named after ``\glsdesc{uobyqa}.''
It aims at solving the nonlinear optimization problem~\eqref{eq:nlc} in the unconstrained case, i.e., when~$\fset = \R^n$.

At iteration $k$, \gls{uobyqa} constructs the model~$\objm$ for the objective function~$\obj$ by the fully determined quadratic interpolation on the interpolation set~$\xpt$ containing~$(n + 1)(n + 2) / 2$ points.
The trust-region subproblem~\eqref{eq:trsp} is solved with the Mor{\'{e}}-Sorensen algorithm~\cite{More_Sorensen_1983}.
For the geometry-improving subproblem~\eqref{eq:biglag}, Powell developed an inexact algorithm that requires only~$\mathcal{O}(n^2)$ operations.
See~\cite[\S~2]{Powell_2002} for more details.

\Gls{uobyqa} updates the trust-region radius~$\rad$ in a noteworthy way.
The update is typical for trust-region methods, except that a lower bound~$\radlb$ is imposed on~$\rad$.
The value of~$\radlb[k]$ can be regarded as an indicator for the current accuracy of the algorithm.
Without imposing~$\rad[k] \ge \radlb[k]$, the trust-region radius~$\rad[k]$ may be reduced to a value that is too small for the current accuracy, making the interpolation points concentrate too much.
The value of~$\radlb[k]$ is never increased and is decreased when the \gls{uobyqa} decides that the work for the current value of~$\radlb[k]$ is finished.
It decides so if~$\rad[k]$ reaches its lower bound~$\radlb[k]$, the current trust-region trial step does not perform well, and the current interpolation set seems adequate for the current accuracy.
See~\cite[\S~3]{Powell_2002} for more information on the updates of~$\rad$ and~$\radlb[k]$.

\subsection{\glsfmttext{newuoa}, \glsfmttext{bobyqa}, and \glsfmttext{lincoa}}
\label{ssec:nbloa}

Later on, based on the underdetermined quadratic interpolation introduced in Subsection~\ref{ssec:iptprob}, Powell developed his last three \gls{dfo} solvers, namely \gls{newuoa}~\cite{Powell_2006,Powell_2008}, \gls{bobyqa}~\cite{Powell_2009}, and \gls{lincoa}.
\Gls{bobyqa} and \gls{lincoa} are named respectively after ``\glsdesc{bobyqa}'' and ``\glsdesc{lincoa},'' but Powell~\cite{Powell_2006,Powell_2008} did not specify the meaning of~\gls{newuoa}, which is likely an acronym for ``\glsdesc{newuoa}.''
It is worth mentioning that Powell \emph{never} published a paper to introduce \gls{lincoa}, and~\cite{Powell_2015} discusses only how to solve its trust-region subproblem.

\Gls{newuoa}, \gls{bobyqa}, and \gls{lincoa} aim at solving unconstrained, bound-constrained, and linearly constrained problems, respectively.
They all set~$\fsetm$ in the trust-region subproblem~\eqref{eq:trsp} to be~$\fset$, corresponding to the whole space for \gls{newuoa}, a box for \gls{bobyqa}, and a polyhedron for \gls{lincoa}.

To solve the trust-region subproblem~\eqref{eq:trsp}, \gls{newuoa} employs the Steihaug-Toint \gls{tcg} algorithm~\cite{Steihaug_1983,Toint_1981}; if the boundary of the trust region is reached, then \gls{newuoa} may make further changes to the trust-region step, each one obtained by searching in the two-dimensional space spanned by the current step and the corresponding gradient of the trust-region model~\cite[\S~5]{Powell_2006}.
\Gls{bobyqa} solves~\eqref{eq:trsp} by an active-set variant of the \gls{tcg} algorithm, and it may also improve the \gls{tcg} step by two-dimensional searches if it reaches the trust-region boundary~\cite[\S~3]{Powell_2009}.
\Gls{lincoa} uses another active-set variant of \gls{tcg} to solve the trust-region subproblem~\eqref{eq:trsp} with linear constraints~\cite[\S~3 and \S~5]{Powell_2015}.
An accessible description of the \gls{tcg} algorithms employed by \gls{bobyqa} and \gls{lincoa} can be found in~\cite[\S\S~6.2.1--6.2.2]{Ragonneau_2022}.
\Gls{newuoa}, \gls{bobyqa}, and \gls{lincoa} manage the trust-region radius in a way similar to \gls{uobyqa}, imposing a lower bound~$\radlb$ on~$\rad$ when updating~$\rad$.

When solving the geometry-improving subproblem~\eqref{eq:biglag}, \gls{newuoa} first takes\break $\iter \pm \radalt (y_k^\drop - \iter)/\norm{y_k^\drop -\iter}$, with the sign that provides the larger value of~$\ell_k^{@\drop}$, and then revises it by a procedure similar to the two-dimensional searches that improve the \gls{tcg} step for~\eqref{eq:trsp} (see~\cite[\S~6]{Powell_2006}).
\Gls{bobyqa} computes two approximate solutions to~\eqref{eq:biglag} and chooses the better one: the first one solves~\eqref{eq:biglag} with an additional constraint that~$\iter[]$ is located on the straight lines through~$\iter$ and another point in~$\xpt$, and the second is obtained by a Cauchy step for~\eqref{eq:biglag} (see~\cite[\S~3]{Powell_2009}).
The geometry-improving step of \gls{lincoa} is more complex, as it is chosen from three approximate solutions to~\eqref{eq:biglag}:
\begin{enumerate}
    \item the point that maximizes~$\abs{\ell_k^{@\drop}}$ within the trust region on the lines through~$\iter$ and another point in~$\xpt$,
    \item a point obtained by a gradient step that maximizes~$\abs{\ell_k^{@\drop}}$ within the trust region, and
    \item a point obtained by a projected gradient step that maximizes~$\abs{\ell_k^{@\drop}}$ within the trust region, the projection being made onto the null space of the constraints that are considered active at~$\iter$.
\end{enumerate}
Note that the first two cases disregard the linear constraints (i.e.~$x\in\fsetm = \fset$), while the third case considers only the active constraints.
\Gls{lincoa} first selects the point among the first two alternatives for a larger value of~$\abs{\ell_k^{@\drop}}$; further, this point is replaced with the third alternative if the latter nearly satisfies the linear constraints~$x\in\fset$ while rendering a value of~$\abs{\ell_k^{@\drop}}$ that is not too small compared with the above one.

\Gls{bobyqa} respects the bound constraints~$x\in \fset$ when solving the trust-region subproblem~\eqref{eq:trsp} and the geometry-improving subproblem~\eqref{eq:biglag}, even though these problems are solved approximately.
It also chooses the initial interpolation set~$\xpt[1]$ within the bounds.
Therefore, \gls{bobyqa} is a feasible method.
In contrast, \gls{lincoa} may violate the linear constraints when solving the geometry-improving subproblem and when setting up the initial interpolation set.
Consequently, \gls{lincoa} is an infeasible method, which requires~$\obj$ to be defined even when the linear constraints are not satisfied.

\section{The \texorpdfstring{\gls{pdfo}}{PDFO} package}
\label{sec:pdfo}

This section details the main features of \gls{pdfo}, in particular the signature of the main function, solver selection, problems preprocessing, bug fixes, and handling failures of function evaluations.
For more features of~\gls{pdfo}, we refer to its homepage at~\url{https://www.pdfo.net}\,.

Before starting, we emphasize that \gls{pdfo} does not re-implement Powell's solvers but rather enables Python and MATLAB to call Powell's Fortran implementation.
At a low level, it uses F2PY\footnote{\url{https://numpy.org/doc/stable/f2py}\,.} to interface Python with Fortran, and MEX to interface MATLAB with Fortran, although users never need such knowledge to employ~\gls{pdfo}.

\subsection{Signature of the main function}

The philosophy of \gls{pdfo} is simple: providing a single function named \texttt{pdfo} to solve \gls{dfo} problems with or without constraints, calling Powell's Fortran solvers in the backend.
It takes for input an optimization problem of the form
\begin{subequations}
    \label{eq:pdfo}
    \begin{align}
        \min_{x \in \R^n}   & ~ \obj(x)\\
        \st                 & ~ \xl \le x \le \xu, \label{eq:pdfo-b}\\
                            & ~ \aub x \le \bub, ~ \aeq x = \beq, \label{eq:pdfo-l}\\
                            & ~ \cub(x) \le 0, ~ \ceq(x) = 0, \label{eq:pdfo-nl}
    \end{align}
\end{subequations}
where~$\obj$ a real-valued objective function, while $\ceq$ and $\cub$ are vector-valued constraint functions.
The bound constraints are given by~$n$-dimensional vectors~$\xl$ and~$\xu$, which may take infinite values.
The linear constraints are formulated by real matrices~$\aeq$ and~$\aub$ together with real vectors~$\beq$ and~$\bub$ of proper sizes.
We allow one or more of the constraints~\eqref{eq:pdfo-b}--\eqref{eq:pdfo-nl} to be absent.
Being a specialization of~\eqref{eq:nlc}, problem~\eqref{eq:pdfo} is broad enough to cover numerous applications of~\gls{dfo}.

In the Python version of~\gls{pdfo}, the signature of the \texttt{pdfo} function is compatible with the \texttt{minimize} function available in the \texttt{scipy.optimize} module of SciPy.
It can be invoked in exactly the same way as \texttt{minimize} except that \texttt{pdfo} does not accept derivative arguments.
The MATLAB version of~\gls{pdfo} designs the \texttt{pdfo} function following the signature of the \texttt{fmincon} function available in the Optimization Toolbox of MATLAB.
In both Python and MATLAB, users can check the detailed syntax of \texttt{pdfo} by the standard \texttt{help} command.

\subsection{Automatic selection of the solver}
\label{subsec:solver-selection}

When invoking the \texttt{pdfo} function, the user may specify which solver to call in the backend.
In the Python and MATLAB versions, this can be done by setting the argument \texttt{method} and the option \texttt{solver}, respectively.
These names are consistent with those used in \texttt{scipy.optimize.minimize} and \texttt{fmincon}.
However, if the user does not specify a solver or chooses a solver that is incapable of solving the problem (e.g., \gls{uobyqa} cannot solve constrained problems), then \texttt{pdfo} selects the solver as follows.
\begin{enumerate}
    \item If the problem is unconstrained, then \gls{uobyqa} is selected when~$2 \le n \le 8$, and \gls{newuoa} is selected when~$n = 1$ or~$n > 8$.
    \item If the problem is bound-constrained, then \gls{bobyqa} is selected.
    \item If the problem is linearly constrained, then \gls{lincoa} is selected.
    \item Otherwise, \gls{cobyla} is selected.
\end{enumerate}
The problem type is detected automatically according to the input.
In the unconstrained case, we select \gls{uobyqa} for small problems because it is more efficient, and the number~\num{8} is set according to our experiments on the CUTEst~\cite{Gould_Orban_Toint_2015} problems.
Note that Powell's implementation of \gls{uobyqa} cannot handle univariate unconstrained problems, for which \gls{newuoa} is chosen.

In addition to the \texttt{pdfo} function, \gls{pdfo} provides functions named \texttt{cobyla}, \texttt{uobyqa}, \texttt{newuoa}, \texttt{bobyqa}, and~\\texttt{lincoa}, which invoke the corresponding solvers directly, but it is highly recommended to call the solvers via the \texttt{pdfo} function.

\subsection{Problem preprocessing}
\label{subsec:pdfo-preprocessing}

\Gls{pdfo} preprocesses the input of the user in order to fit the data structure expected by Powell's Fortran code.

For example, \gls{lincoa} needs a feasible starting point to work properly unless the problem is infeasible.
If the starting point is not feasible, then \gls{lincoa} would modify the right-hand sides of the linear constraints to make it feasible and then solve the modified problem.
Therefore, for linearly constrained problems, \gls{pdfo} attempts to project the user-provided starting point onto the feasible region before passing the problem to the Fortran code so that a feasible problem will not be modified by~\gls{lincoa}.

Another noticeable preprocessing of the constraints made by \gls{pdfo} is the treatment of the linear equality constraints in~\eqref{eq:pdfo-l}.
As long as these constraints are consistent, we eliminate them and reduce~\eqref{eq:pdfo} to an~$(n - \rank \aeq)$-dimensional problem.
This is done using a QR factorization of~$\aeq$.
The main motivation for this reduction comes again from~\gls{lincoa}, which accepts only linear inequality constraints.
An alternative approach is to write a linear equality constraint as two inequalities, but our approach reduces the dimension of the problem, which is beneficial for the efficiency of \gls{dfo} solvers in general.

\subsection{Bug fixes in the Fortran source code}
\label{subsec:bug-corrections}

The current version of \gls{pdfo} patches several bugs in the original Fortran source code, particularly the following ones.
\begin{enumerate}
    \item The solvers may encounter infinite loops.
        This happens when the exit conditions of a loop can never be met because variables involved in these conditions become NaN due to floating point exceptions.
        The user's program will never end if this occurs.
    \item The Fortran code may encounter memory errors due to uninitialized indices.
        This is because some indices are initialized according to conditions that can never be met due to NaN, similar to the previous case.
        The user's program will crash if this occurs.
\end{enumerate}

In our extensive tests based on the CUTEst problems, these bugs take effect from time to time but not often.
They are activated only when the problem is rather ill-conditioned or the inputs are rather extreme.
This has been observed, for instance, on the CUTEst problems DANWOODLS, GAUSS1LS, and LAKES with some perturbation and randomization.

Even though these bugs are rarely observed in our tests, it is vital to patch them for two reasons.
First, their consequences are severe once they occur.
Second, application problems are often more irregular and savage than the testing problems we use, and hence the bugs may be triggered more often than we expect.
Nevertheless, \gls{pdfo} allows the users to call Powell's original Fortran code without these patches by setting the option \texttt{classical} to true, which is highly discouraged.

\subsection{Handling failures of function evaluations}
\label{ssec:barrier}

\Gls{pdfo} tolerates NaN values returned by function evaluations.
Such a value can be used to indicate failures of function evaluations, which are common in applications of~\gls{dfo}.

To cope with NaN values, \gls{pdfo} applies a \emph{moderated extreme barrier}.
Suppose that~$\obj(\tilde{x})$ is evaluated to NaN at a certain~$\tilde{x}\in\R^n$.
\Gls{pdfo} takes the view that~$\tilde{x}$ violates a hidden constraint~\cite{LeDigabel_Wild_2023,Audet_Caporossi_Jacquet_2020}.
Hence it replaces NaN with a large but finite number \texttt{HUGEFUN} (e.g., $10^{30}$) before passing~$\obj(\tilde{x})$ to the Fortran solver, so that the solver can continue to progress while penalizing~$\tilde{x}$.
Indeed, since Powell's solvers construct trust-region models by interpolation, all points that are close to~$\tilde{x}$ will be penalized.
Similar things are done when the constraint functions return NaN.
A caveat is that setting~$\obj(\tilde{x})$ to \texttt{HUGEFUN} may lead to extreme values or even NaN in the coefficients of the interpolation models, but Powell's solvers turn out to be quite tolerant of such values.

The original extreme barrier approach~\cite[Equation~(13.2)]{Conn_Scheinberg_Vicente_2009b} sets \texttt{HUGEFUN} to~$\infty$, which is inappropriate for methods based on interpolation.
In fact, we also moderate~$\obj(\tilde{x})$ to~$\texttt{HUGEFUN}$ if it is actually evaluated to~$\infty$.
Our approach is clearly naive, but it is better than terminating the solver once the function evaluation fails.
In our experiments, this simple approach significantly improves the robustness of~\gls{pdfo} with respect to failures of function evaluation, as will be demonstrated in Subsection~\ref{ssec:nan}.
There do exist other more sophisticated approaches~\cite{Audet_Caporossi_Jacquet_2020}, which will be explored in the future.

\section{Numerical results}
\label{sec:numerical}

This section presents numerical experiments on \gls{pdfo}.
Since Powell's solvers are widely used as benchmarks in \gls{dfo}, extensive comparisons with standard \gls{dfo} solvers are already available in the literature~\cite{More_Wild_2009,Rios_Sahinidis_2013}.
Instead of repeating such comparisons, the purpose of our experiments is the following.
\begin{enumerate}
    \item Demonstrate the fact the \gls{pdfo} is capable of adapting to noise without fine-tuning according to the noise, in contrast to methods based on finite differences.
        This is done in Subsection~\ref{ssec:noise} by comparing \gls{pdfo} with finite-difference \gls{cg} and \gls{bfgs} on unconstrained CUTEst problems.
    \item Verify the effectiveness of the moderated extreme barrier mentioned in Subsection~\ref{ssec:barrier} for handling failures of function evaluations.
        This is done in Subsection~\ref{ssec:nan} by testing \gls{pdfo} with and without the barrier on unconstrained CUTEst problems.
    \item Illustrate the potential of \gls{pdfo} in hyperparameter optimization problems from machine learning, echoing the observations made in~\cite{Ghanbari_Scheinberg_2017} about trust-region \gls{dfo} methods for such problems.
        This is done in Subsection~\ref{ssec:hypertune} by comparing~\gls{pdfo} with two solvers from the Python package \texttt{hyperopt}.\footnote{\url{https://hyperopt.github.io/hyperopt}\,.}
\end{enumerate}

Our experiments are carried out in double precision based on the Python version of \gls{pdfo}~2.2.0.
The finite-difference \gls{cg} and \gls{bfgs} are provided by SciPy~1.11.3.
The version of \texttt{hyperopt} is~0.2.7.
All these packages are tested with the latest stable version at the time of writing.
We conduct the test on a ThinkStation P620 with an AMD Ryzen Threadripper PRO 3975WX CPU and 64 GB of memory, the operating system being Ubuntu~22.04, and the Python version being~3.10.12.

\subsection{Stability under noise}
\label{ssec:noise}

We first compare \gls{pdfo} with finite-difference \gls{cg} and \gls{bfgs} on unconstrained problems with multiplicative Gaussian noise.
We take the view that multiplicative noise makes more sense if the scale of the objective function changes widely, as is often the case in applications.

SciPy provides both~\gls{cg} and \gls{bfgs} under the \texttt{minimize} function in the \texttt{scipy.optimize} module.
These methods rely on the gradient of the objective function.
If we do not provide such information, SciPy approximates it by finite differences, and we do the same in our experiments.
For \gls{pdfo}, we specify \gls{newuoa} as the solver, while setting all the other options to the default ones.
In particular, the initial trust-region radius is~$1$, the final trust-region radius is~$10^{-6}$, and the number of interpolation points is~$2n + 1$, where~$n$ is the dimension of the problem being solved.
We perform the comparison on~\num{166} unconstrained problems with~$n \le 50$ from the CUTEst~\cite{Gould_Orban_Toint_2015} problem set using PyCUTEst~1.5.1~\cite{Fowkes_Roberts_Burmen_2022}.
For each testing problem, the starting point is set to the one provided by CUTEst, and the maximal number of function evaluations is~$500n$.

Let~$\sigma \ge 0$ be the noise level to test.
For a testing problem with the objective function~$\obj$, we define
\begin{equation}
    \label{eq:noisy-obj}
    \tilde{\obj}_\sigma(\iter[]) = [1 + \sigma R(\iter[])] \obj(\iter[]),
\end{equation}
with $R(\iter[])\sim \NN(0, 1)$ being independent and identically distributed when~$x$ varies.
If~$\sigma = 0$, then~$\tilde{f}_\sigma = f$, corresponding to the noise-free case.
In general,~$\sigma$ is the standard deviation of the noise.
The function~$\tilde{\obj}_\sigma$ is the one received by the optimization solvers.
In particular, the \gls{cg} and \gls{bfgs} methods of SciPy approximate the components of~$\nabla \tilde{\obj}_\sigma(\iter[])$ by the finite difference
\begin{equation*}
    \frac{\partial \tilde{\obj}_\sigma}{\partial [{\iter[]}]_i}(\iter[]) \approx \frac{\tilde{\obj}_\sigma(\iter[] + h e_i) - \tilde{\obj}_\sigma(\iter[])}{h},
\end{equation*}
where~$[{\iter[]}]_i$ denotes the~$i$th component of~$\iter[]$,~$e_i \in \R^n$ is the~$i$th standard coordinate vector, and~$h > 0$ is the difference parameter
\begin{equation}
    \label{eq:scipy-h}
    h = \sgn([{\iter[]}]_i) \max \set{\abs{[{\iter[]}]_i}, 1} \sqrt{u},
\end{equation}
where~$u \approx 2.2 \times 10^{-16}$ is the unit roundoff.
When there is noise ($\sigma > 0$), we also test another value of~$h$, namely
\begin{equation}
    \label{eq:adaptive-h}
    h = \sqrt{\sigma \max \set[\big]{\abs[\big]{\tilde{\obj}_\sigma(\iter[])}, 1}}.
\end{equation}
This adaptive difference parameter is inspired by the optimal choice that depends on the second-order information of~$\obj$ (see, e.g.,~\cite{More_Wild_2012} and~\cite[Equation~(2.2)]{Shi_Etal_2022a}).
However, it relies on the knowledge of the noise level~$\sigma$.
In contrast, \gls{pdfo} does not require the knowledge of~$\sigma$ and as we will see, provides better performance.

Given a noise level~$\sigma\ge 0$ and a convergence tolerance~$\tau\in(0,1)$, we will plot the performance profiles~\cite{More_Wild_2009} of the solvers on the testing problems.
We run all the solvers on all the problems, every objective function being evaluated by its contaminated version~\eqref{eq:noisy-obj}.
For each solver, the performance profile displays the proportion of problems solved with respect to the normalized cost to solve the problem up to the convergence tolerance~$\tau$.
For each problem, the cost to solve the problem is the number of function evaluations needed to achieve
\begin{equation}
    \label{eq:cvt}
    \obj(\iter[0]) - \obj(\iter) \ge (1 - \tau) [\obj(\iter[0]) - \obj_{\ast}],
\end{equation}
and the normalized cost is this number divided by the minimum cost of this problem among all solvers; we define the normalized cost as infinity if the solver fails to achieve~\eqref{eq:cvt} on this problem.
Here,~$\iter[0]$ represents the starting point, and~\cite[\S~2.1]{More_Wild_2009} suggests that the value~$\obj_{\ast}$ should be the least value of~$\obj$ obtained by all solvers.
Note that the convergence test~\eqref{eq:cvt} uses the values of~$\obj$ and not those of~$\tilde{\obj}_\sigma$.
This means that we assess the solvers according to the true objective function values, even though the objective function fed to the solvers is~$\tilde{\obj}_\sigma$, which is contaminated unless~$\sigma = 0$.

To make our results more reliable, when~$\sigma>0$, the final performance profile is obtained by averaging the profiles obtained via the above procedure over ten independent runs.
In addition, the value~$\obj_{\ast}$ in the convergence test~\eqref{eq:cvt} is set to the least value of~$\obj$ obtained by all solvers during all these ten runs plus a run with~$\sigma = 0$.
Finally, for better scaling of the profiles, we plot the binary logarithm of the normalized cost on the horizontal axis, instead of the normalized cost itself.

\begin{figure}[!htbp]
    \begin{subfigure}{.48\textwidth}
        \centering
        \includegraphics[width=\textwidth]{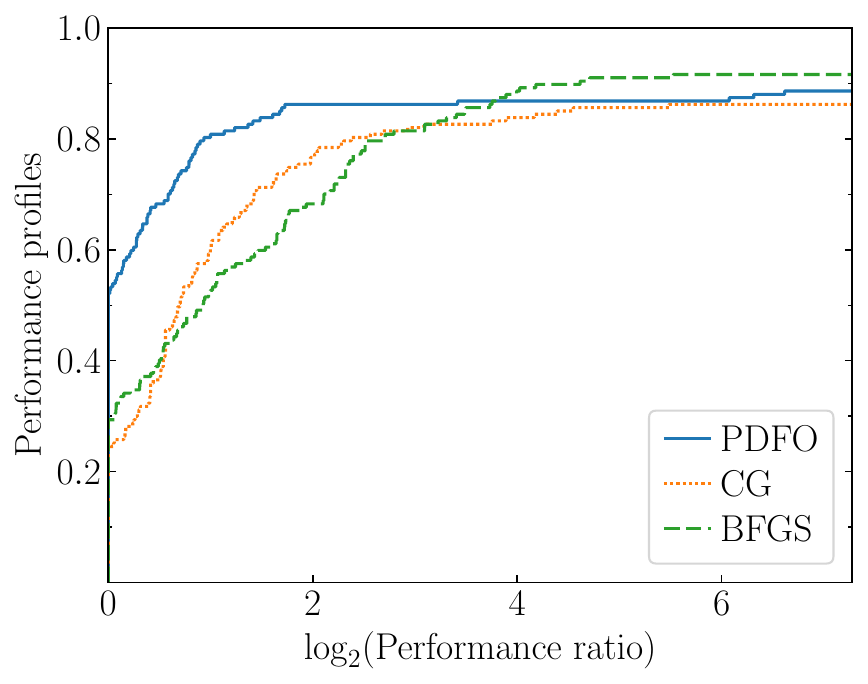}
        \caption{$\sigma = 0$, $\tau = 10^{-2}$}
    \end{subfigure}
    \hfill
    \begin{subfigure}{.48\textwidth}
        \centering
        \includegraphics[width=\textwidth]{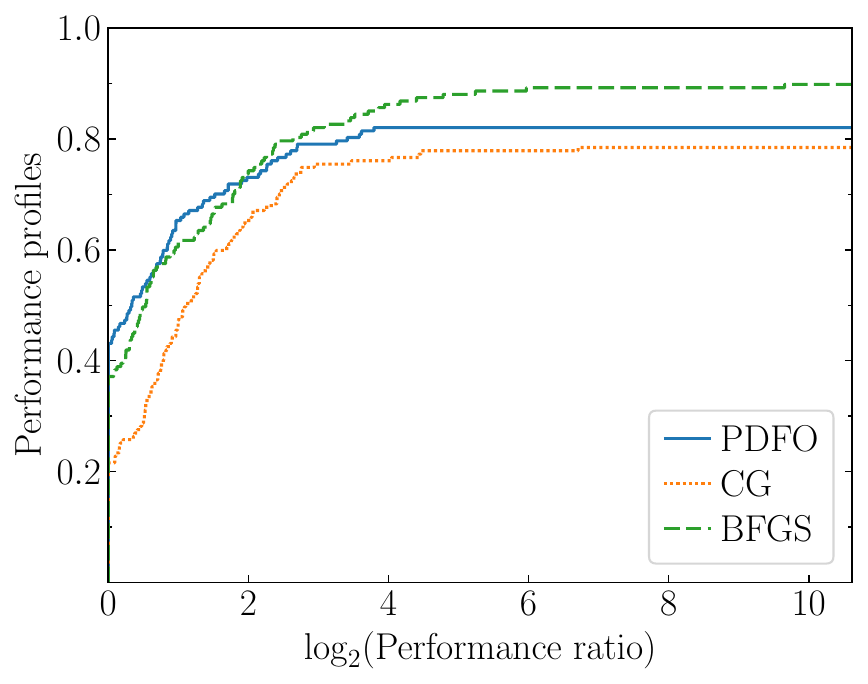}
        \caption{$\sigma = 0$, $\tau = 10^{-4}$}
    \end{subfigure}
    \hfill
    \begin{subfigure}{.48\textwidth}
        \centering
        \includegraphics[width=\textwidth]{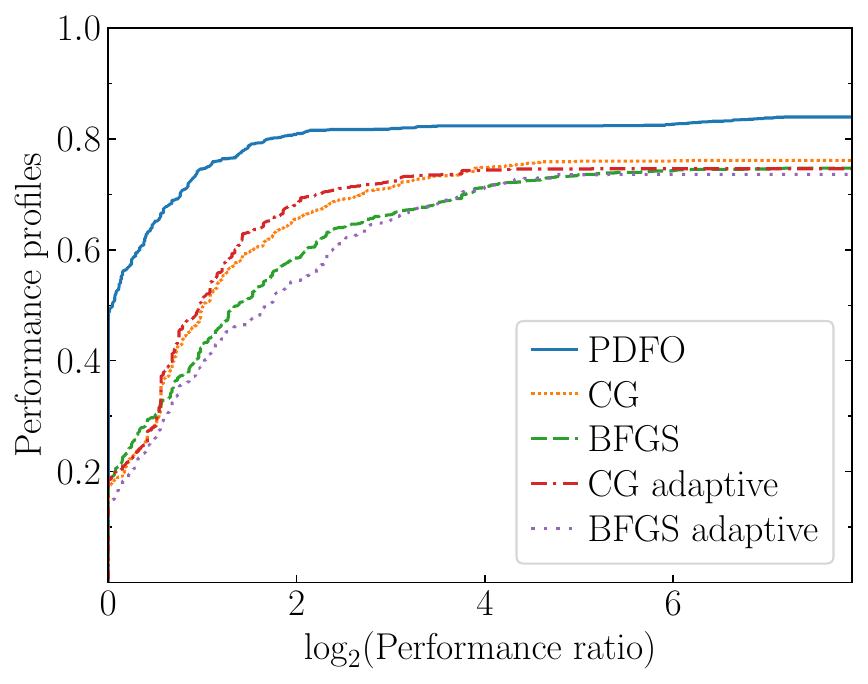}
        \caption{$\sigma = 10^{-10}$, $\tau = 10^{-2}$}
    \end{subfigure}
    \hfill
    \begin{subfigure}{.48\textwidth}
        \centering
        \includegraphics[width=\textwidth]{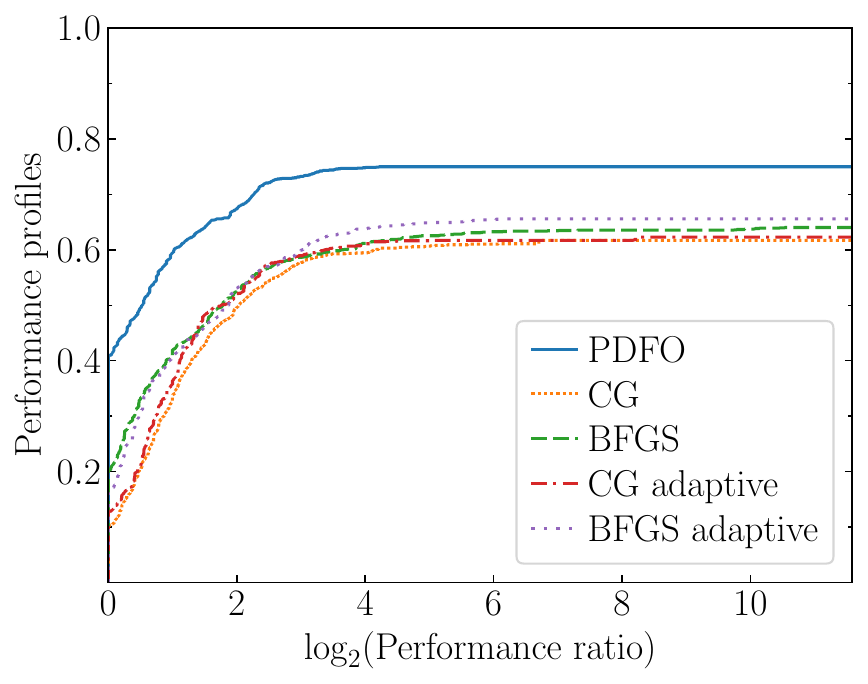}
        \caption{$\sigma = 10^{-10}$, $\tau = 10^{-4}$}
    \end{subfigure}
    \hfill
    \begin{subfigure}{.48\textwidth}
        \centering
        \includegraphics[width=\textwidth]{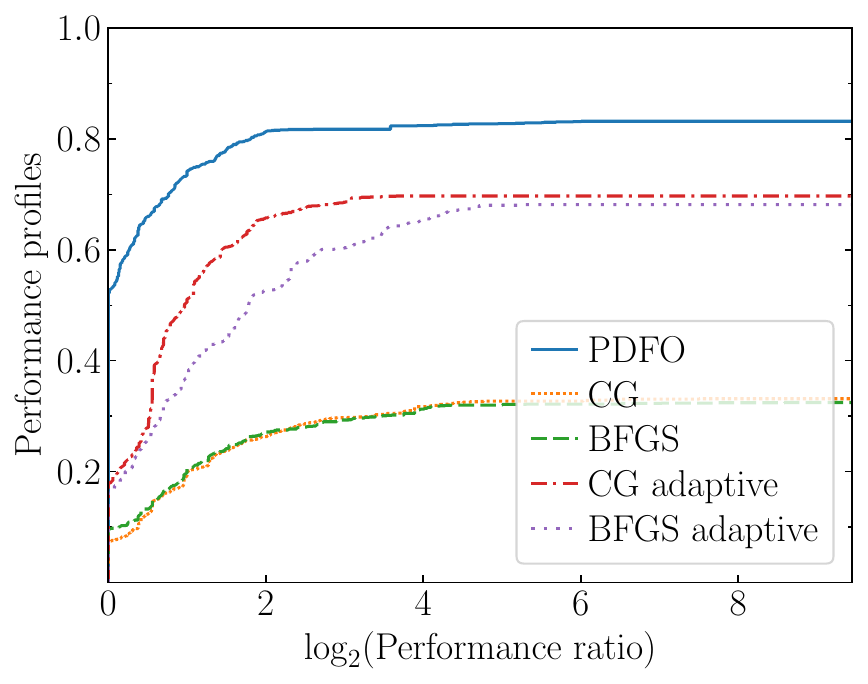}
        \caption{$\sigma = 10^{-8}$, $\tau = 10^{-2}$}
    \end{subfigure}
    \hfill
    \begin{subfigure}{.48\textwidth}
        \centering
        \includegraphics[width=\textwidth]{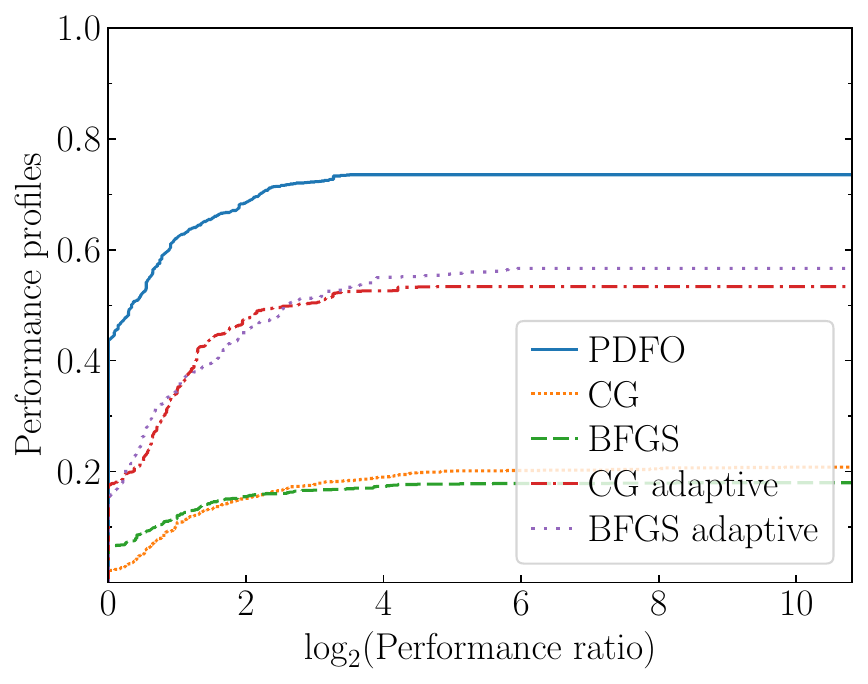}
        \caption{$\sigma = 10^{-8}$, $\tau = 10^{-4}$}
    \end{subfigure}
    \caption{Performance profiles of \gls{pdfo}, \gls{cg}, and \gls{bfgs} on unconstrained problems with the objective functions evaluated by~$\tilde{\obj}_\sigma$ in~\eqref{eq:noisy-obj}}
    \label{fig:noise}
\end{figure}

Figure~\ref{fig:noise} shows the performance profiles of the solvers for the noise levels~$\sigma = 0$, $\sigma = 10^{-10}$, and~$\sigma = 10^{-8}$.
Two profiles are included for each noise level, with the convergence tolerance being~$\tau = 10^{-2}$ and~$\tau = 10^{-4}$ respectively.
In the legends, ``\gls{cg}'' and ``\gls{bfgs}'' denote the finite-difference methods using~\eqref{eq:scipy-h}, whereas ``\gls{cg} adaptive'' and ``\gls{bfgs} adaptive'' are their counterparts with~$h$ given by~\eqref{eq:adaptive-h}.
Note that the subfigures corresponding to~$\sigma = 0$ do not include ``\gls{cg} adaptive'' and ``\gls{bfgs} adaptive'' because these methods are designed to tackle noisy problems.

In the noise-free case ($\sigma = 0$), \gls{pdfo} is more efficient than finite-difference \gls{cg} and \gls{bfgs}, although the distinction is less visible when~$\tau$ is smaller, and~\gls{bfgs} can solve slightly more problems than~\gls{pdfo}.
When there is noise ($\sigma >0$), the advantage of \gls{pdfo} becomes significant.
The performances of \gls{cg} and \gls{bfgs} using~\eqref{eq:scipy-h} deteriorate considerably under noise, even though the noise level is not high and the convergence tolerance is not demanding.
As expected, the adaptive difference parameter~\eqref{eq:adaptive-h} improves the performance of \gls{cg} and \gls{bfgs} by a large margin, but they are still inferior to \gls{pdfo}.
This shows the advantage of \gls{pdfo} in noisy settings, even when compared to methods that explicitly use the magnitude of the noise.
We have conducted similar experiments with larger values of~$\sigma$.
The results are similar to those shown in Figure~\ref{fig:noise}, except that the advantage of \gls{pdfo} over \gls{cg} and \gls{bfgs} with~\eqref{eq:adaptive-h} is even more visible, and the other two methods barely solve any problem.
For conciseness, we do not include these results.

It is not surprising that \gls{cg} and \gls{bfgs} perform unfavorably when using~\eqref{eq:scipy-h}.
As mentioned in~\cite{More_Wild_2012,Shi_Etal_2022a,Shi_Etal_2023}, the difference parameter should be adapted according to the noise level, as we do in the adaptive versions.
In contrast, \gls{pdfo} does not need such adaptation when handling noisy problems.

To summarize, the performance of finite-difference \gls{cg} and \gls{bfgs} is encouraging when there is no noise, yet much more care is needed when the problems are noisy.
In contrast, \gls{pdfo} adapts to noise automatically in our experiment, demonstrating good stability under noise without requiring knowledge about the noise level.
This is because Powell's methods (\gls{newuoa} in this experiment) gradually adjust the geometry of the interpolation set during the iterations, making progress until the interpolation points are too close to distinguish noise from true objective function values.
This is not specific to Powell's methods but also applies to other algorithms that sample the objective function on a set of points with adaptively controlled geometry, including finite-difference methods with well-chosen difference parameters~\cite{Shi_Etal_2022a}.

\subsection{Robustness with respect to failures of function evaluations}
\label{ssec:nan}

We now test the robustness of the solvers when function evaluations fail from time to time.
We assume that the objective function returns NaN if the evaluation fails, which occurs randomly with a certain probability.
As mentioned in Section~\ref{ssec:barrier}, \gls{pdfo} uses a moderated extreme barrier to handle such failures.
To verify the effectiveness of this approach, we compare \gls{pdfo} with its variant that does not apply the barrier.
To make the experiment more informative, we also include the finite-difference \gls{cg} and~\gls{bfgs} tested before, which do not handle evaluation failures particularly.
The solvers are set up in the same way as in the previous experiment, and we still employ the~\num{166} unconstrained CUTEst problems used previously.

Let~$p \in [0,1]$ be the failure probability of function evaluations.
For a testing problem with the objective function~$\obj$, we define
\begin{equation}
    \label{eq:nan-obj}
    \hat{\obj}_p(x) = \begin{cases}
        \obj(x)     & \text{if~$U(x) \ge p$},\\[0.5ex]
        \text{NaN}  & \text{otherwise},
    \end{cases}
\end{equation}
where $U(\iter[])$ follows the uniform distribution on~$[0,1]$, being independent and identically distributed when~$x$ varies.
Note that~$\hat{\obj}_0 = f$.
In the experiment, the solvers can evaluate~$\obj$ only via~$\hat{\obj}_p$.
We plot the performance profiles of the solvers in a way that is similar to the previous experiment.
The profiles are also averaged over ten independent runs.
For each problem, the value~$\obj_{\ast}$ in the convergence test~\eqref{eq:cvt} is set to the least value of~$\obj$ obtained by all solvers during these ten runs plus a run with~$p = 0$.

\begin{figure}[!htbp]
    \begin{subfigure}{.48\textwidth}
        \centering
        \includegraphics[width=\textwidth]{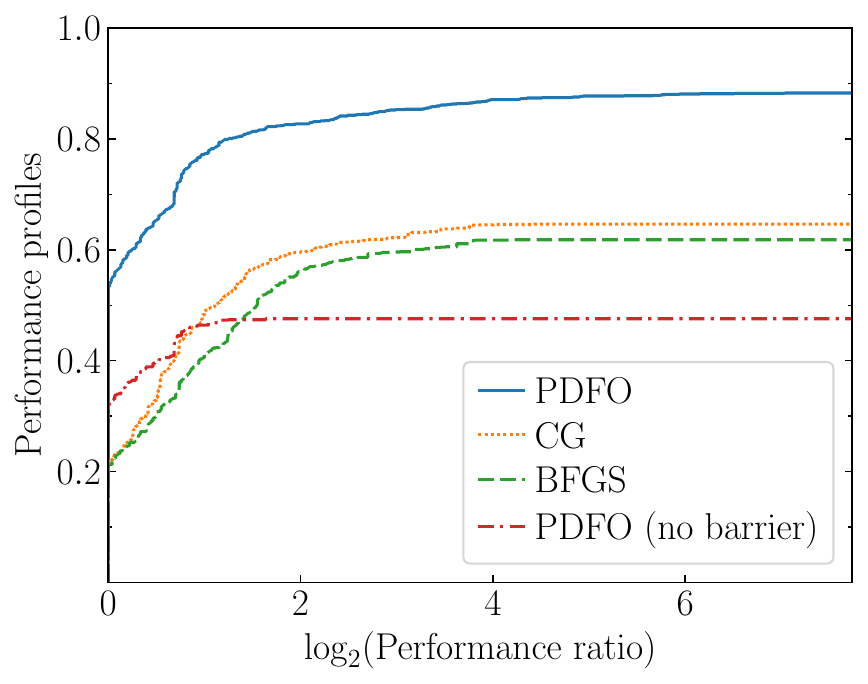}
        \caption{$p = 0.01$, $\tau = 10^{-2}$}
    \end{subfigure}
    \hfill
    \begin{subfigure}{.48\textwidth}
        \centering
        \includegraphics[width=\textwidth]{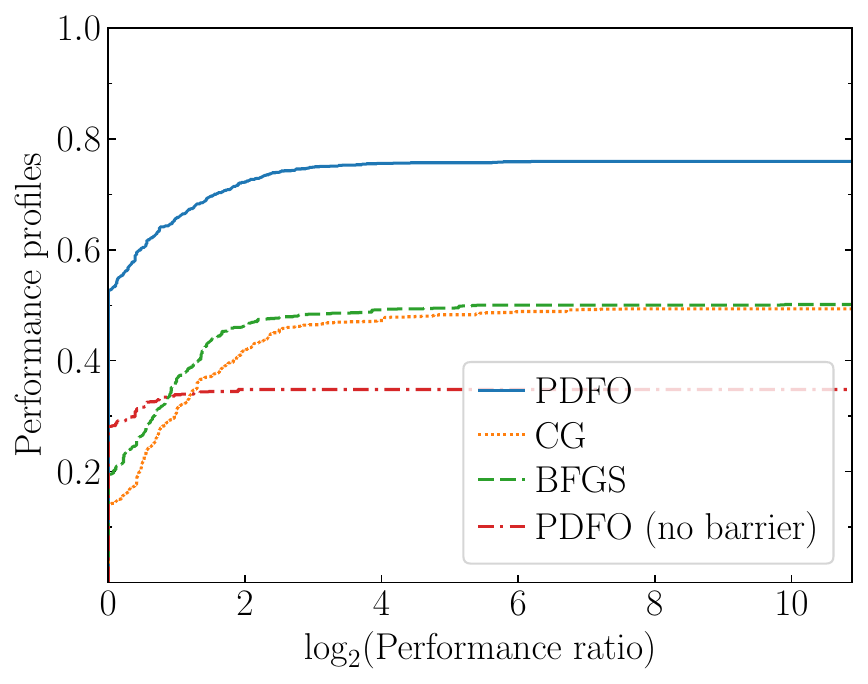}
        \caption{$p = 0.01$, $\tau = 10^{-4}$}
    \end{subfigure}
    \hfill
    \begin{subfigure}{.48\textwidth}
        \centering
        \includegraphics[width=\textwidth]{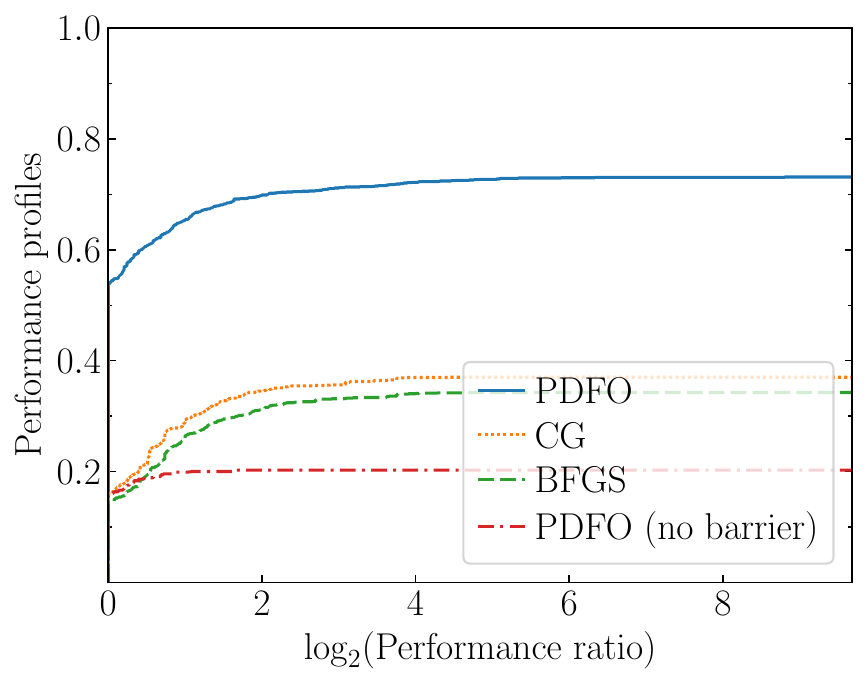}
        \caption{$p = 0.05$, $\tau = 10^{-2}$}
    \end{subfigure}
    \hfill
    \begin{subfigure}{.48\textwidth}
        \centering
        \includegraphics[width=\textwidth]{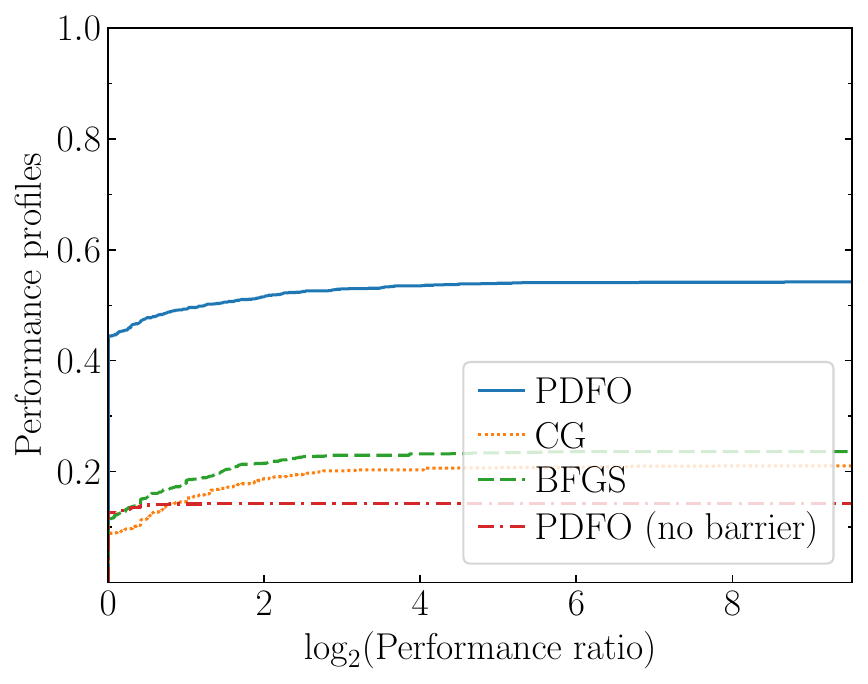}
        \caption{$p = 0.05$, $\tau = 10^{-4}$}
    \end{subfigure}
    \caption{Performance profiles of \gls{pdfo}, \gls{cg}, \gls{bfgs}, and~\gls{pdfo} without barrier on unconstrained problems with the objective functions evaluated by~$\hat{\obj}_p$ in~\eqref{eq:nan-obj}}
    \label{fig:nan}
\end{figure}

Figure~\ref{fig:nan} shows the performance profiles of the solvers with~$p = 0.01$ and~$p=0.05$.
Two profiles are included for each~$p$, with the convergence tolerance being~$\tau = 10^{-2}$ and~$\tau = 10^{-4}$ respectively.

The contrast is clear.
Compared with finite-difference \gls{cg} and \gls{bfgs}, \gls{pdfo} is more efficient and solves significantly more problems given the same convergence tolerance.
Moreover, comparing \gls{pdfo} and its no-barrier counterpart, we see that the moderated extreme barrier improves evidently the robustness of~\gls{pdfo} with respect to failures of function evaluations, even though it is a quite naive approach.
When~$p = 0.05$, the function evaluation fails roughly once every $20$ times, but \gls{pdfo} can still solve almost~$55\%$ of the problems up to the convergence tolerance~$\tau = 10^{-4}$ in the sense of~\eqref{eq:cvt}, whereas all its competitors solve less than $25\%$.
We speculate that the moderated extreme barrier will also benefit other model-based \gls{dfo} methods, including those based on finite differences.
It deserves further investigation in the future.

\subsection{An illustration of hyperparameter optimization with~\glsfmttext{pdfo}}
\label{ssec:hypertune}

We now consider a hyperparameter optimization problem from machine learning and illustrate the potential of \gls{pdfo} for such problems.
We compare \gls{pdfo} with \gls{rs} and \gls{tpe}, two solvers from the Python package \texttt{hyperopt} for hyperparameter optimization.

Our experiment is inspired by~\cite[\S~5.3]{Ghanbari_Scheinberg_2017}, which investigates the application of trust-region \gls{dfo} methods to hyperparameter optimization.
Similar to~\cite[\S~5.3]{Ghanbari_Scheinberg_2017}, we tune the $C$-SVC model detailed in~\cite[\S~2.1]{Chang_Lin_2011} for binary classifications.
This model relies on two hyperparameters: a penalty parameter~$C \in (0, \infty)$ and a kernel parameter~$\gamma \in (0, \infty)$.
As suggested by~\cite[\S~9]{Chang_Lin_2011}, we tune~$C$ and~$\gamma$ for the performance of the~$C$-SVC.
We model this process as solving the problem
\begin{equation}
    \label{eq:hypertune}
    \max ~ P(C, \gamma) \quad \st \quad C > 0, ~ \gamma > 0,
\end{equation}
where~$P(C, \gamma)$ measures the performance corresponding to parameters~$(C, \gamma)$.
In our experiment, we define~$P$ based on the AUC score~\cite[\S~3]{Ghanbari_Scheinberg_2017}, which lies in~$[0,1]$ and measures the quality of a classifier on a dataset, the higher the better.
More precisely, $P(C, \gamma)$ is set to a five-fold cross-validation AUC score as follows.
Split the training dataset~$\mathcal{S}$ into five folds, and train the $C$-SVC five times, each time on a union of four distinct folds.
After each training, calculate the AUC score of the resulting classifier on the fold not involved in the training, leading to five scores, the average of which is~$P(C, \gamma)$.

Our experiment is based on binary classification problems from LIBSVM,\footnote{\url{https://www.csie.ntu.edu.tw/~cjlin/libsvmtools/datasets}\,.} where we adopt three datasets detailed in Table~\ref{tab:htdata}.
LIBSVM divides each dataset~$\mathcal{D}$ into two disjoint subsets, namely a training dataset~$\mathcal{S}$ and a testing dataset~$\mathcal{T}$.
The training in the evaluation of~$P$ is done using the \texttt{SVC} class of the Python package scikit-learn.\footnote{\url{https://scikit-learn.org/stable/modules/generated/sklearn.svm.SVC.html}\,.}
We solve~\eqref{eq:hypertune} by~\gls{pdfo}, \gls{rs}, and~\gls{tpe} to obtain the tuned parameters~$(\bar{C}, \bar{\gamma})$.
As in~\cite[\S~5.3]{Ghanbari_Scheinberg_2017}, we modify the constraints of~\eqref{eq:hypertune} to~$C\in[10^{-6}, 1]$ and~$\gamma\in [1, 10^{3}]$.
For better scaling of the problem, we perform the maximization with respect to~$(\log_{10}C,\, \log_{10}\gamma)$ instead of~$(C, \gamma)$, the initial guess being chosen randomly from~$[-6, 0]\times[0, 3]$.
The solver of~\gls{pdfo} is \gls{bobyqa}, for which we set the maximal number of function evaluations to~\num{100}.
For~\gls{rs} and~\gls{tpe}, we try both~\num{100} and~\num{300} for the maximal number of function evaluations, and they do not terminate until this number is reached.

\begin{table}[!htb]
    \caption{Datasets from LIBSVM}
    \label{tab:htdata}
    \centering
    \begin{tabular}{lcS[table-format=5.0]S[table-format=5.0]}
        \toprule
        \multicolumn{1}{c}{Dataset} & {Number of features}  & {Size of~$\mathcal{S}$}   & {Size of~$\mathcal{T}$}\\
        \midrule
        splice                      & 60                    & 1000                      & 2175\\
        svmguide1                   & 4                     & 3088                      & 4000\\
        ijcnn1                      & 22                    & 49990                     & 91701\\
        \bottomrule
    \end{tabular}
\end{table}

To assess the quality of the tuned parameters~$(\bar{C}, \bar{\gamma})$, we train our model on~$\mathcal{S}$ with\break $(C, \gamma)=(\bar{C}, \bar{\gamma})$, and calculate both the AUC score and accuracy of the resulting classifier on~$\mathcal{T}$, the latter being the fraction of correctly classified data points.
Note that~$\mathcal{T}$ is not involved in the tuning process.
Table~\ref{tab:ht} presents the results for this experiment, where \#$P$ denotes the number of evaluations of the function~$P$ and ``Time'' is the computing time for obtaining~$(\bar{C}, \bar{\gamma})$.

\begin{table}[!htb]
    \caption{Hyperparameter tuning using \gls{pdfo}, \gls{rs}, and \gls{tpe}}
    \label{tab:ht}
    \centering
    \begin{tabular}{llS[table-format=1.3]S[table-format=1.3]cS[scientific-notation=true,round-mode=places,table-format=1.2e1]}
        \toprule
        \multicolumn{1}{c}{Dataset} & \multicolumn{1}{c}{Solver}    & {AUC Score} & {Accuracy}    & {\#$P$}   & {Time~(\si{\second})}\\
        \midrule
        \multirow{5}{*}{splice}     & \gls{pdfo}                    & 0.927       & 0.737         & 33        & 6.56\\
        \cmidrule{2-6}
                                    & \multirow{2}{*}{\gls{rs}}     & 0.500       & 0.520         & 100       & 9.16\\
                                    &                               & 0.500       & 0.520         & 300       & 27.4\\
        \cmidrule{2-6}
                                    & \multirow{2}{*}{\gls{tpe}}    & 0.500       & 0.520         & 100       & 8.89\\
                                    &                               & 0.500       & 0.520         & 300       & 26.4\\
        \midrule
        \multirow{5}{*}{svmguide1}  & \gls{pdfo}                    & 0.995       & 0.966         & 51        & 3.45\\
        \cmidrule{2-6}
                                    & \multirow{2}{*}{\gls{rs}}     & 0.994       & 0.961         & 100       & 11.3\\
                                    &                               & 0.995       & 0.968         & 300       & 34.0\\
        \cmidrule{2-6}
                                    & \multirow{2}{*}{\gls{tpe}}    & 0.995       & 0.965         & 100       & 9.53\\
                                    &                               & 0.995       & 0.968         & 300       & 25.3\\
        \midrule
        \multirow{5}{*}{ijcnn1}     & \gls{pdfo}                    & 0.997       & 0.980         & 44        & 2440\\
        \cmidrule{2-6}
                                    & \multirow{2}{*}{\gls{rs}}     & 0.997       & 0.982         & 100       & 3939\\
                                    &                               & 0.997       & 0.976         & 300       & 11366\\
        \cmidrule{2-6}
                                    & \multirow{2}{*}{\gls{tpe}}    & 0.998       & 0.979         & 100       & 3387\\
                                    &                               & 0.998       & 0.979         & 300       & 8361\\
        \bottomrule
    \end{tabular}
\end{table}

In terms of the AUC score and accuracy, \gls{pdfo} achieves a clearly better result than \gls{rs} and \gls{tpe} on the ``splice'' dataset, and they all attain comparable results on the other datasets.
However, \gls{pdfo} always uses much fewer function evaluations, and hence, much less computing time.
The difference in the computing time is particularly visible on the ``ijcnn1'' dataset, as each evaluation of~$P$ takes much time due to the large data size.

Note that our intention is not to manifest that \gls{pdfo} outperforms existing approaches for hyperparameter tuning in general, which is unlikely the case.
Indeed, \gls{pdfo} has limitations in handling hyperparameter tuning problems, as many such problems contain discrete variables and cannot be handled by \gls{pdfo} directly.
Our objective is rather to provide an example that shows the possibility of applying Powell's methods to hyperparameter optimization, which is not well studied up to now.
In doing so, we also hope to call for more investigation on \gls{dfo} methods for machine learning problems in general, as is suggested in~\cite{Ghanbari_Scheinberg_2017}.

\section{Concluding remarks}
\label{sec:conclude}

We have presented the \gls{pdfo} package, which aims at simplifying the use of Powell's \gls{dfo} solvers by providing user-friendly interfaces.
More information about the package can be found on the homepage of the package at \url{https://www.pdfo.net}\,, including the detailed syntax of the interfaces, extensive documentation of the options, and several examples to illustrate the usage.

In addition, we have provided an overview of Powell's methods behind \gls{pdfo}.
The overview does not intend to repeat Powell's description of the methods, but rather to provide a summary of the main features and structures of the methods, highlighting the intrinsic connections and similarities among them.
We hope that the overview will ease the understanding of Powell's methods, in the same way as the \gls{pdfo} package eases the use of these methods.

Besides Powell's solvers, \gls{pdfo} also provides a unified interface for \gls{dfo} solvers.
Such an interface can facilitate the development and comparison of different \gls{dfo} solvers.
The interface can readily accommodate solvers other than those by Powell, for example, the \gls{cobyqa} (\glsdesc{cobyqa}) solver for general nonlinearly constrained \gls{dfo} problems (see~\cite[Chapters~5--7]{Ragonneau_2022} and~\cite{Ragonneau_Zhang_cobyqa}).

Finally, we stress that \gls{pdfo} does not implement Powell's \gls{dfo} solvers in Python or MATLAB, but only interfaces Powell's implementation with such languages.
The implementation of these solvers in Python, MATLAB, and other languages is a project in progress under the name of \gls{prima} (\glsdesc{prima})~\cite{Zhang_prima}.

\paragraph{Acknowledgments.}

This paper corresponds to Chapter~3 of the PhD thesis of Tom M.\break Ragonneau~\cite{Ragonneau_2022}, co-supervised by Zaikun Zhang and Professor Xiaojun Chen from The Hong Kong Polytechnic University.
Both authors are very grateful to Professor Chen for her support, encouragement, and guidance during the PhD studies.
Also, the authors would like to thank Professor Ya-xiang Yuan for his everlasting encouragement and support.
Finally, the authors thank the editors and referees for their suggestions, which have substantially improved the software and manuscript.

\printbibliography

\end{document}